\documentclass [a4paper,twoside,12pt]{article}

\usepackage[utf8]{inputenc}
\usepackage{amsmath,amsfonts,amssymb}
\usepackage{vmargin,graphicx,theorem}
\usepackage[english]{babel}
\usepackage{enumerate}
\usepackage{color}

\RequirePackage[colorlinks,linkcolor=blue,citecolor=blue,urlcolor=blue]{hyperref}
\usepackage{amssymb, amsmath, amsfonts, latexsym}
\usepackage{enumerate,graphicx}
\usepackage{pst-fill,pst-grad,pst-plot,pst-eucl,pstricks-add,pst-node}

\setpapersize[portrait]{A4}
\setmarginsrb{1.5cm}{1cm}{1.5cm}{2.5cm}{1.5cm}{0cm}{0.5cm}{2cm}

\def\I{1\!{\rm l}}
\def\tr{\textmd{trace}\,}

\def\Ric{ {\rm{Ric}}} 

\def\Id{{\rm{Id}}} 

\newcommand{\bbR}{{\mathbb {R}}}
\newcommand{\bbS}{{\mathbb {S}}} 

\newcommand{\disp}{\displaystyle}

\newcommand{\dR}{\ensuremath{\mathbb{R}}}

\newtheorem{theo}{Theorem}

\newcommand{\proofend}{~$\rhd$}
\newcommand{\proofbegin}{~$\lhd$}

\newcommand{\PAR}[1]{\ensuremath{{\left(#1\right)}}} 
\newcommand{\SBRA}[1]{\ensuremath{{\left[#1\right]}}} 

\renewcommand{\phi}{\varphi}

\renewcommand{\geq}{\geqslant}

\def\tr{\mathop{\rm tr}\nolimits} 

\newcommand{\upchi}{\raise1pt\hbox{$\chi$}}

\newcommand{\be}{\begin{equation}}
\newcommand{\ee}{\end{equation}}
\def\benu{\begin{enumerate}}
\def\eenu{\end{enumerate}}

\newcommand{\bthm}{\begin{theo}}
\newcommand{\ethm}{\end{theo}}

\newcommand{\bcor}{\begin{corollaire}}
\newcommand{\ecor}{\end{corollaire}}
\newcommand{\bdefi}{\begin{definition}}
\newcommand{\edefi}{\end{definition}}
\newcommand{\bprop}{\begin{proposition}}
\newcommand{\eprop}{\end{proposition}}

\newcommand{\blem}{\begin{lemme}}
\newcommand{\elem}{\end{lemme}}

\newcommand{\beqna}{\begin{eqnarray}}
\newcommand{\eeqna}{\end{eqnarray}}

\newcommand{\beqnas}{\begin{eqnarray*}}
\newcommand{\eeqnas}{\end{eqnarray*}}

\newtheorem{lemme}[theo]{Lemma}
\newtheorem{definition}[theo]{Definition}
\newtheorem{proposition}[theo]{Proposition}
\newtheorem{corollaire}[theo]{Corollary}

\newcommand{\bpf}{\begin{preuve}}
\newcommand{\epf}{ \end{preuve} \medskip}
\newenvironment{preuve}{\noindent{\it Proof. --- }}
{\hfill\rule{1.3mm}{2mm}\par} 
\newcommand{\bcas}{\begin{cases}}
\newcommand{\ecas}{\end{cases}}
\newcommand{\brmq}{\begin{rmq}}
\newcommand{\ermq}{\end{rmq}}
\newcommand{\brmqs}{\begin{rmqs}}
\newcommand{\ermqs}{\end{rmqs}}
\newtheorem{rmq}[theo]{Remark}
\newtheorem{rmqs}[theo]{Remarks}

\newcommand{\beq}{\begin{equation}}\newcommand{\eeq}{\end{equation}}
\title{Sharp Beckner-type inequalities for Cauchy and spherical distributions}
\author{Dominique Bakry\thanks{Institut de Math\'ematiques de Toulouse, Umr Cnrs 5219, Universit\'e de Toulouse - Paul Sabatier, bakry@math.univ-toulouse.fr}, Ivan Gentil\thanks{Univ Lyon, Université Claude Bernard Lyon 1, CNRS UMR 5208, Institut Camille Jordan, 43 blvd. du 11 novembre 1918, F-69622 Villeurbanne cedex, France. gentil@math.univ-lyon1.fr} and Gr\'egory Scheffer\thanks{SCOR Investment Partners 5, avenue Kléber, 75795 Paris cedex 16, gregory.scheffer@gmail.com}}

\date{\today}
\begin{document}
\maketitle

\abstract{Using some harmonic extensions on the upper-half plane, and probabilistic representations, and curvature-dimension inequalities with some negative dimensions,  we obtain some new optimal functional inequalities of the Beckner type  for the Cauchy type distributions on the Euclidean space.  These optimal inequalities appear to be equivalent to some non tight optimal Beckner
inequalities on the sphere, and the family appears to be a new form of the Sobolev inequality.}
\bigskip

\medskip


\section{Introduction}
The so-called $\Gamma_2$ criterium is a way to prove functional inequalities such as Poincar\'e or logarithmic Sobolev inequalities.  For instance,  let $\Psi:\bbR^d\mapsto \bbR$ be a smooth function such that $\nabla^2\psi\geq \rho\,{\rm I}$ with $\rho>0$, then the probability measure 
$$
d\mu_\psi=\frac{e^{-\psi}}{Z}dx,
$$
where $Z$ is the normalization constant which turns $\mu_\psi$ into a probability, satisfies both the Poincar\'e inequality,  
$$
\int f^2d\mu_\psi-\PAR{\int fd\mu_\psi}^2\leq \frac{1}{\rho}\int |\nabla f|^2d\mu_\psi,
$$
and the logarithmic Sobolev  inequality
$$
\int f^2\log f^2d\mu_\psi-\int f^2d\mu_\psi\log \int f^2d\mu_\psi\leq \frac{2}{\rho}\int |\nabla f|^2d\mu_\psi,
$$
for any smooth function $f$. More generally,   the same hypothesis leads to  a family of interpolation inequalities between Poincar\'e and logarithmic Sobolev, namely,    for any $p\in (1,2]$, 
$$
\frac{p}{p-1}\PAR{\int f^2d\mu_\psi-\PAR{\int |f|^{2/p}d\mu_\psi}^p}\leq \frac{2}{\rho}\int |\nabla f|^2d\mu_\psi,
$$
where the case $p=2$ corresponds to the Poincar\'e inequality and the limit  $p\rightarrow 1$ to the logarithmic Sobolev inequality.  One of the technics used to prove such inequalities are based on the $\Gamma_2$ calculus (the $CD(\rho,+\infty)$ condition) introduced by the first author with M. \'Emery in~\cite{bakryemery} (see also~\cite{bgl-book}). This family is the so-called  Beckner inequalities, proved for the usual Gaussian measure,  that is when $\psi= \frac{\|x\|^2}{2}$, by W. Beckner in~\cite{beckner1} (in that case  $\rho=1$). Some improvements of such inequalities can be found in~\cite{arnolddolbeault05,dns,bolley-gentil}. Let use note also that recent developments for the Beckner inequalities applied to the Gaussian distribution in an optimal way have been proved in~\cite{17-volberg},  also for uniformly strictly log-concave measures. Let us finish to quote the recent article of Nguyen~\cite{18-nguyen} where the author generalizes some results of this paper with a completely different method, the method is extended later in~\cite{18dgz}.

On the other hand, on the sphere $\bbS^d\subset\bbR^d$, with  similar arguments, the following Beckner inequalities can be stated, $p\in(1,2]$,  
\begin{equation}
\label{eq-last}
\frac{p}{p-1}\PAR{\int f^2d\mu_\bbS-\PAR{\int |f|^{2/p}d\mu_\bbS}^p}\leq \frac{2}{d}\int |\nabla f|^2d\mu_\bbS,
\end{equation}
see~\cite{bgl-book},  where in this situation, $\rho= d-1$.

It is important to notice that all Beckner inequalities (for the spherical and the Gaussian models)  are optimal in the sense that constants in front of the right hand side are optimal. However,  only the constant functions saturate these inequalities when $p\neq 2$.

\medskip

The aim of this article is two-fold. First, we explain and improve in a general context the method introduced by the last author~\cite{scheffer} to obtain some new optimal inequalities. This method is mainly probabilistic, and makes strong use of various processes associated with the underlying structure behind the inequalities.  He used this method to obtain an optimal Poincar\'e inequality under a $CD(\rho,n)$ condition, which is a stronger form of the condition $\nabla^2\psi \geq \rho Id$, and without integration by parts.  More precisely, the method is enough robust to be applied to a non-symmetric operator even if our example is the usual Laplacian which is symmetric.    This method is rather technical and  is based on the $CD(\rho,n)$ condition with a negative dimension $n$, together with a construction of some sub-harmonic functional.  Let us note that negative dimension on the curvature-dimension $CD(\rho,n)$ condition has been introduced and used in~\cite{scheffer,ohta,milman}.  Stochastic calculus is also a main tool of our approach.

Secondly, improving unpublished results of the last author,  we apply this method to obtain Beckner inequalities  for the  generalized Cauchy distribution, that is the  probability measure on $\bbR^d$ with density 
$$
d\nu_b(x)=\frac{1}{Z(1+|x|^2)^{b}}dx
$$
where $b>d/2$ and $Z$ is the normalization constant. This measure does not satisfy any $CD(\rho,+\infty)$ condition with $\rho\geq0$  but some $CD(0,n)$ 
condition with $n<0$, as will be explained with in detail in the paper. We prove the following new  optimal family of inequalities,
$$
\frac{p}{p-1}\left[\int {f^2}d\nu_b-\left(\int {|f|^{2/p}}d\nu_b\right)^{p}\right]\leq
 \frac{1}{(b-1)}\int \Gamma(f)(1+|y|^2)d\nu_b.
$$
where $b\geq d+1$ and $p\in[1+1/(b-d),2]$.

Of course, when $p=2$ we recover a weighted Poincar\'e inequality for the Cauchy distribution proved in many papers, see~\cite{scheffer-phd,blanchet2007,bonforte2010,nguyen-14}.

But more surprisingly, this family is equivalent to the following one on the sphere $\bbS^d$,  
\begin{equation}
\label{last?}
\int h^{2}d\mu_\bbS\leq A\left(\int |h|^{2/p}d\mu_\bbS\right)^{p}
+\frac{16}{(m+2-d)(3d-2+m)}\int \Gamma_\bbS(h)d\mu_\bbS,
\end{equation}
where
$$m\geq d+2, \quad
p=1+\frac{2}{m-d}\in[1,2]
$$
and $A\geq1$ is some explicit constant to be described below, converging to 1 as $p$ converges to 1, or equivalently when $m$ converges to infinity.  Inequality~\eqref{eq-last} is a tight inequality, that is constant functions are optimal, on the other hand inequality~\eqref{last?} is a non-tight inequality. On the other hand, this inequality~\eqref{last?} is optimal with explicit extremal functions. Moreover this new family of inequalities, or more precisely its behavior when $p$ converges to 1,  appears as a new form of a Sobolev inequality on the sphere, as many other ones described for example in~\cite{bcls}. 

\medskip

The paper is organized as follows. In the next section, we describe the curvature-dimension condition $CD(\rho,n)$, even for negative $n$, and derive a  general form of associated sub-harmonic functionals. In Section~\ref{sec-Qtm}, we define an operator related to the generalized Cauchy distribution, and show that it  satisfies the generalized curvature-dimension condition. In Section~\ref{sec-be-cauchy}, we prove the new Beckner inequalities for these generalized Cauchy distribution. Finally, in Section~\ref{sec-be-sphere}, we prove the non-tight Becker inequalities on the $d$-dimensional sphere and show how they relate to the Sobolev inequality.

\bigskip

{\bf Notations:} In all this paper, $d$ will be the dimension of the main space, and satisfies $d\geq1$, and $d\geq2$ in sections~\ref{sec-tightness} and \ref{sec-be-sphere}. For every $X,Y\in\bbR^d$, $X\cdot Y=\sum_{i=1}^dX_iY_i$ is the usual scalar product and $|X|^2$ is the Euclidean norm in $\bbR^d$.

\section{General properties on curvature-dimension condition}

Let $(M^d,\mathfrak{g})$ be a smooth connected $d$-dimensional  Riemannian manifold and $\Delta_{\mathfrak{g}}$ its Laplace-Beltrami operator associated. The associated Ricci tensor is denoted $\Ric_\mathfrak{g}$. For any diffusion operator $L$ on $M$ (that is a second order semi-elliptic differential operator with no zero-order term, see~\cite{bgl-book}) we define the so-called  carr\'e du champ   operator 
$$
\Gamma^L(f,g)=\frac{1}{2}\SBRA{L(fg)-fLg-gLf}
$$ 
with $\Gamma^L(f,f)=\Gamma^L(f)$ and its iterated operator 
$$
\Gamma_{2}^L(f,f)=\Gamma_{2}^L(f)=\frac{1}{2}L(\Gamma_L(f))-\Gamma^L(f,Lf),
$$ 
 for any smooth functions $f$ and $g$. We note $\Gamma$ (resp.  $\Gamma_2$) instead $\Gamma^L$ (resp. $\Gamma_2^L$) when there is no possible confusion.
\subsection{Definitions}
 \begin{definition}[$CD(\rho,n)$ condition]
An operator $L=\Delta_{\mathfrak{g}}+X$ with $X$ a smooth vector field satisfies a $CD(\rho,n)$ condition with $\rho\in\bbR$ and $n\in\bbR\setminus[0,d)$ if 
\begin{equation}
\label{eq-def-gamma2}
\Gamma_2(f)\geq \rho \Gamma(f)+ \frac{1}{n} (Lf)^2,
\end{equation}
for any smooth function $f$. 
\end{definition}
As we shall see of the proof in Lemma~\ref{lem-saintflour}, this condition will never hold when $n\in[0,d)$ but we will extend it for $n=0$.

\begin{lemme}[\cite{bakrystflour}]
\label{lem-saintflour}
For any $\rho\in\bbR$ and $n\notin[0,d]$ the operator $L=\Delta_{\mathfrak{g}}+X$ satisfies a $CD(\rho,n)$ if and only if 
\beq
\label{eq-condition-tensor}
\frac{n-d}{n} (\Ric(L)-\rho \mathfrak{g})\geq \frac{X\otimes X}{n},
\eeq
and when $n=d$ it reduces only to $X=0$ and the condition becomes $\Ric_\mathfrak{g} \geq \rho \mathfrak{g}$. In the definition, $\Ric(L)=\Ric_\mathfrak{g}-\nabla_S X$ where $\nabla_S X$ is the symmetrized tensor of  $\nabla X$. 
\end{lemme}
\bpf
For completeness, we give a sketch of proof of this result which can be helpful for the rest of the paper. For any smooth function $f$, the Bochner-Lichnerowicz-Weitzenbock formula states that (cf.~\cite[Sec C.5]{bgl-book}),
\beq
\label{eq-bochner}
\Gamma_2(f)\!=\!||\nabla^2f||_{H.S.}^2+\Ric_\mathfrak{g}(\nabla f,\nabla f)-\nabla X_S(\nabla f,\nabla f)=||\nabla^2f||_{H.S.}^2+\Ric(L)(\nabla f,\nabla f),
\eeq
where $||\nabla^2f||_{H.S.}$ is the Hilbert-Schmidt norm of $f$. On a fixed point $x\in M$, we choose a local chart such that the 
metric $\mathfrak{g}=(\mathfrak{g}_{i,j})$ is  the identity 
matrix at the point $x$. In that case, $\nabla^2f=(\nabla_{i,j}f)$ is a $d$-dimensional symmetric matrix denoted $Y$, $\Gamma(f)=\sum_{i=1}^d(\partial_if)^2$
 and let denote $Z=\nabla f$.  Since $\sum_{i,j=1}^dY_{i,j}^2=\tr(Y^2)$, the $CD(\rho,n)$ condition~\eqref{eq-def-gamma2} becomes  
$$
\tr(Y^2)+Z\cdot RZ\geq \rho Z\cdot Z+\frac{1}{n}(\tr(Y)+X\cdot Z)^2,
$$ 
where $R$ is a matrix representing the coordinates of the tensor $\Ric(L)$ at $x$. One can choose a new base  in $\bbR^d$ such that $Y$ is diagonal. If we denote again
 $Z$ in the new base in $\bbR^d$, we have to prove such inequality, for every vectors $(Y_{ii})_{1\leq i\leq d}$ and $(Z_i)_{1\leq i\leq d}$,  
$$
\sum_{i=1}^dY_{ii}^2+Z\cdot RZ\geq \rho Z\cdot Z+\frac{1}{n}(\sum_{i=1}^dY_{ii}+X\cdot Z)^2.
$$ 
By taking $Y_{ii}$ all equals, from $\sum_{i=1}^dY_{ii}^2\geq\frac{1}{d}(\sum_{i=1}^dY_{ii})^2$, the previous inequality  is equivalent to 
$$
\frac{1}{d}(\sum_{i=1}^dY_{ii})^2+Z\cdot RZ\geq \rho Z\cdot Z+\frac{1}{n}(\sum_{i=1}^dY_{ii}+X\cdot Z)^2.
$$ 
Let $y=\sum_{i=1}^dY_{ii}$, then the inequality becomes 
$$
\frac{y^2}{d}+Z\cdot RZ\geq \rho Z\cdot Z+\frac{1}{n}(y+X\cdot Z)^2.
$$ 
In other words, the $CD(\rho,n)$ condition is satisfied if and only if, for every $y\in\bbR$ and, $Z\in\bbR^d$
$$
y^2\PAR{\frac{1}{d}-\frac{1}{n}}-y\frac{2}{n}X\cdot Z+Z\cdot (R-\rho\text{Id})Z-\frac{1}{n}(X\cdot Z)^2\geq0.
$$ 

Then, either $n=d$, then we need to assume that $X=0$ and $Z\cdot (R-\rho\text{Id})Z\geq0$,
either $\frac{1}{d}-\frac{1}{n}>0$ (with translates to $n\notin[0,d]$)  and for any vector $Z$, 
$$
\frac{(X\cdot Z)^2}{n}\leq \frac{n-d}{n}Z\cdot (R-\rho\text{Id})Z,
$$
that is condition~\eqref{eq-condition-tensor}. 
\epf
One may know define the tensor curvature-dimension condition, $TCD(\rho,n)$.   
\begin{definition}[$TCD(\rho,n)$ condition and quasi-models $QM(\rho,n)$]
Let $\rho\in\bbR$, $n\in\bbR\setminus\{d\}$ and an operator $L=\Delta_{\mathfrak{g}}+X$ with $X$ a smooth vector field. 
\begin{itemize}
\item $L$ satisfies a $TCD(\rho,n)$ condition if 
\begin{equation}
\label{eq-def-tcd}
\frac{n-d}{n}\big(\Ric(L)-\rho \mathfrak{g}\big)\geq \frac{X\otimes X}{n},
\end{equation}
when $n\neq0$,  
and 
$$
-d(\Ric(L)-\rho \mathfrak{g})\leq X\otimes X,
$$
when $n=0$. 
\item The operator $L$ is a quasi model $QM(\rho,n)$ if 
\beq
\label{eq-quasi-models}
(n-d)(\Ric(L)-\rho \mathfrak{g})={X\otimes X}.
\eeq
\end{itemize}
\end{definition}
Indeed, this new definition allows us to extend~\eqref{eq-def-gamma2} when $n\in[0,d)$, and we shall see later crucial examples of such operators, which are even quasi-models.

\subsection{Sub-harmonic functionals}
 In this section, we construct sub-harmonic functionals from harmonic functions under $TCD(\rho, n)$. We systematically explore the most general ones, in view of further use. 
 
 The next result is rather technical, and we shall see in the next sections how it leads to Poincaré
 inequalities or Beckner inequalities.

\bthm[Sub-harmonic functionals under $TCD(\rho ,n)$ or $QM(\rho,n)$]
\label{grad.bound.gal.Phi} 
~

Let $\Phi:I\times[0,\infty)\mapsto \bbR$ be a smooth function ($I$ open interval of $\bbR$), and denote $\Phi_i$ and $\Phi_{ij}$  its first and second derivatives, $i,j\in\{1,2\}$. Let $L= \Delta_{\mathfrak{g}} +X$, where $X$ is a smooth vector field. 

\medskip

Let $\rho\in\bbR$ and $n\in\bbR\setminus\{d\}$. We assume that $L$ satisfies the  $TCD(\rho ,n)$ condition with $n\in\bbR\setminus(0,d]$ (or the 
$QM(\rho,n)$ condition). Then, if  $F:M\mapsto I$ is smooth and satisfies $L(F)=0$, then $L\SBRA{\Phi (F,\Gamma(F))} \geq 0$  
as soon as $\Phi$ satisfies at any point $(y,z)\in I\times[0,\infty)$, either  
\beq
\label{eq-condition-gradient}
\bcas 
\disp\Phi_2(y,z) >0; \\
\disp\Phi_2(y,z)\frac{n+1-d}{n-d}+ 2z\Phi_{22}(y,z) >0;\\
\disp (n-d)\big(n\Phi _{2}(y,z)+2(n-1)z\Phi _{22}(y,z)\big)> 0;\\
\disp 2\rho \Phi_2(y,z)+ \Phi_{11}(y,z)- 2(n-1) \frac{z\Phi_{12}^2(y,z) }{n\Phi_2(y,z)+ 2 (n-1)z\Phi_{22}(y,z) } \geq 0,
\ecas
\eeq
either
\beq
\label{eq-55}
\bcas 
\disp\Phi_2(y,z)=0;\\
\disp z\Phi_{22}(y,z)=0; \\
\disp z\Phi_{12}(y,z)=0;\\
z\Phi_{11}(y,z)\geq0,\\
\ecas
\eeq
either 
\beq
\label{eq-55}
\bcas 
\disp\Phi_2(y,z)=0;\\
z\Phi_{22}(y,z)>0;\\
\Phi_{11}(y,z)\Phi_{22}(y,z)-\Phi_{12}^2(y,z)\geq0.\\
\ecas
\eeq
\ethm

\brmqs
\begin{itemize}
\item For simplicity, we omit the variables $y$ and $z$ inside the function $\Phi$. When $n\notin [0,d]$, condition $\Phi_2\frac{n+1-d}{n-d}+ 2\Phi_{22} z >0$ is contained into 
$\Phi _{2}+2\frac{n-1}{n}z\Phi _{22}> 0$. Hence, provided that $n\notin [0,d]$, then~\eqref{eq-condition-gradient} boils down to
\beq
\label{eq-condition-0d}
\bcas
 \Phi_2 >0, \\
\disp\Phi _{2}+2\frac{n-1}{n}z\Phi _{22}> 0,\\
\disp 2\rho \Phi_2+ \Phi_{11}- 2(n-1) \frac{\Phi_{12}^2z }{n\Phi_2+ 2 (n-1)\Phi_{22} z} \geq 0.
\ecas
\eeq 
It is interesting to notice that when $n\notin [0,d]$ the dimension $d$ does not appear in the previous conditions~\eqref{eq-condition-0d}. 
\item As pointing out by a referee, condition~\eqref{eq-condition-gradient} has a nice formulation by using a new variable coming from~\cite{18Ivanisvili-volberg}. The authors set  $M(y,u)=\Phi(x,\sqrt{u})$, then condition~\eqref{eq-condition-gradient} is equivalent to, for any $(u,y)\in\bbR\times (0,\infty)$,
$$
(n-d)\left(
\begin{array}{cc}
\disp M_{11}(y,u)+\frac{\rho}{u}M_{12}(y,u) & M_{12}(y,u)\\
\disp M_{12}(y,u) & M_{22}(y,u)+\frac{1}{(n-1)u}M_{2}(y,u)
\end{array}
\right)
\geq0,
$$
together with the conditions $M_2(y,u)>0$ and $M_{22}(y,u)+\frac{1}{(n-1)u}M_{2}(y,u)>0$ where we use again the notation $M_{ij}=\partial^2_{ij}M$. Actually, we do not know implications of this new formulation of our condition~\eqref{eq-condition-gradient}.
\end{itemize}
\ermqs
\bpf 
The proof is an improvement of the proof of Lemma~\ref{lem-saintflour}. 

If $L F=0$, then from the diffusion property of $L$ (cf.~\cite[Sec~1.14]{bgl-book}),  
\begin{multline}
L(\Phi(F, \Gamma(F)))=\Phi_2 L\Gamma(F) + \Phi_{11} \Gamma(F)+ 2\Phi_{12} \Gamma(F, \Gamma(F))+ \Phi_{22} \Gamma(\Gamma(F), \Gamma(F))\\ \label{eq-lgamma}
=2\Phi_2 \Gamma_2(F) + \Phi_{11} \Gamma(F)+ 2\Phi_{12} \Gamma(F, \Gamma(F))+ \Phi_{22} \Gamma(\Gamma(F), \Gamma(F)),
\end{multline}
since $L\Gamma(F)= 2\Gamma_2(F)$ when $LF=0$,  for simplicity we note $\Phi_{ij}=\Phi_{ij}(F,\Gamma(F))$.
At some point $x\in M$, we write all the expressions appearing in $L(\Phi(F, \Gamma(F)))$. From~\eqref{eq-bochner}, they make use only of $\nabla\nabla F$ and $\nabla F$. We chose a system of coordinates in which at this point the metric $\mathfrak{g}$ is identity,  $\nabla\nabla F$ is diagonal, with eigenvalues $(\lambda_i)_{1\leq i\leq d}$. The components of $\nabla F$ are denoted $Z_i$, and the components of $X$ are denoted $-X_i$.

Since $L(F)=0$,  the constraint can be written,
$$
\sum_{i=1}^d \lambda_i = \sum_{i=1}^d X_iZ_i.
$$
and 
$$
\Gamma(F, \Gamma(F))= 2\sum_{i=1}^d \lambda_i Z_i^2, \quad~\Gamma(\Gamma(F), \Gamma(F))= 4\sum_{i=1}^d \lambda_i^2 Z_i^2,\quad \Gamma(F)=\sum_{i=1}^d Z_i^2.
$$
From the $TCD(\rho,n)$ conditions and~\eqref{eq-bochner}, the $\Gamma_2$ operator has a  lower bound,  
$$
\Gamma_2(F)= \sum_{i=1}^d \lambda_i^2+ \Ric(L)(Z,Z)\geq \sum_{i=1}^d \lambda_i^2+ \rho\sum_{i=1}^d Z_i^2+ \frac{1}{n-d}\Big(\sum_{i=1}^d X_iZ_i\Big)^2.
$$
Under the $TCD(\rho,n)$ this is an inequality and we need to assume at this stage that $n\notin(0,d]$ and when the operator is a quasi-model $QM(\rho,n)$ we only need to assume that $n\neq d$ and this is an equality.

Since $\Phi_2\geq 0$, from~\eqref{eq-lgamma}, we just have to  minimize the following expression,  
\begin{multline*}
\frac{L(\Phi(F, \Gamma(F)))}{2}\geq2\Phi_2 \Big(\sum_{i=1}^d \lambda_i^2+ \rho\sum_{i=1}^d Z_i^2+ \frac{1}{n-d}\Big(\sum_{i=1}^d X_iZ_i\Big)^2\Big)\\
 + \Phi_{11}\sum_{i=1}^d Z_i^2+ 4\Phi_{12} \sum_{i=1}^d \lambda_i Z_i^2+ 4\Phi_{22}\sum_{i=1}^d \lambda_i^2 Z_i^2,
\end{multline*}
for all $(\lambda_i)$ and $Z_i$, under the constraint  $\sum_{i=1}^d \lambda_i= \sum_{i=1}^d X_iZ_i$.

Let $z=\sum_{i=1}^d Z_i^2=\Gamma(F)$ and $Z_i^2 = z_i$, then the expression becomes  by using the constraint (after dividing by $2$)
$$
\Phi_2 \Big(\sum_{i=1}^d \lambda_i^2+ \rho z+ \frac{1}{n-d}\Big(\sum_{i=1}^d \lambda_i\Big)^2\Big) + \frac{\Phi_{11}}{2}z+ 2\Phi_{12} \sum_{i=1}^d \lambda_i z_i+ 2\Phi_{22}\sum_{i=1}^d \lambda_i^2 z_i.
$$

We first optimise on the simplex $\{ z_i\geq 0, \sum_{i=1}^d z_i= z\}$ where $z>0$ is fixed. Of course, we have to minimize along the simplex and the constraint $\{\sum_{i=1}^d\pm X_i\sqrt{z_i}=\sum_{i=1}^d\lambda_i\}$ which enlarge the set. We only minimize along the simplex which is lower. In fact it does not change the result and at this stage we can forget about the constraint.  We first observe that the expression is affine in $z_i$. The vector  $
(z_i)_{1\leq i\leq d}$ belongs to the simplex $\{z_i\geq 0, \sum_{i=1}^d z_i= z\}$, so that the extremals of this expression ($(\lambda_i)_{1\leq i\leq d} $ and 
$z>0$ are fixed) is obtained at the one of the extreme point of the simplex. That is one of the $z_i$ is $z$ and all the other  are $0$. Assuming that 
$z_1=z\neq 0$,  one is now led to minimize
\begin{multline*}
\lambda_1^2 (\Phi_2+ 2 \Phi_{22} z) + 2 \Phi_{12} \lambda_1 z + \Phi_2\sum_{i=2}^d \lambda_i^2 + \frac{\Phi_2}{n-d} \Big(\sum_{i=1}^d \lambda_i\Big)^2+  
z\Big( \Phi_2\rho  + \frac{\Phi_{11}}{2}\Big)\\
=\lambda_1^2 (\Phi_2\frac{n+1-d}{n-d}+ 2 \Phi_{22} z) +  \lambda_1\Big(2 \Phi_{12} z +2\frac{\Phi_2}{n-d} \sum_{i=2}^d \lambda_i\Big)\\
+\Phi_2\Big( \sum_{i=2}^d \lambda_i^2 +\frac{1}{n-d} \Big(\sum_{i=2}^d \lambda_i\Big)^2\Big) +  
z\Big( \Phi_2\rho  + \frac{\Phi_{11}}{2}\Big)
\end{multline*}
Then, we minimize in $\lambda_1$. This imposes the second condition of~\eqref{eq-condition-gradient}, that is  
$$
\Phi_2\frac{n+1-d}{n-d}+ 2\Phi_{22} z >0,
$$
recall that $z=\Gamma(F)$. And the minimizer is given by 
$$
\frac{L(\Phi(F, \Gamma(F)))}{2}\geq\Phi_2\Big(\sum_{i=2}^d \lambda_i^2 + \frac{1}{n-d} \Big(\sum_{i=2}^d \lambda_i\Big)^2 \Big)+z\Big( \Phi_2\rho  + \frac{\Phi_{11}}{2}\Big)- \frac{ (\Phi_{12}z + \frac{\Phi_2}{n-d}\sum_{i=2}^d \lambda_i)^2}{\Phi_2 \frac{n+1-d}{n-d}+ 2 \Phi_{22}z} 
$$
Let assume that $\Phi_2>0$. Decomposing the vector $\Lambda=(\lambda_i)_{2\leq i\leq d}$ into a vector parallel to $(1, \cdots, 1)$ and a vector $u=(u_i)_{2\leq i\leq d}$ orthogonal to it, that is $\Lambda=\lambda(1,\cdots,1)+u$ with $\sum_{i=2}^du_i=0$.  From this decomposition, we have $\sum_{i=2}^d \lambda_i=(d-1)\lambda$ and 
$$
\sum_{i=2}^d \lambda_i^2=\sum_{i=2}^d (\lambda+u_i)^2\geq (d-1)\lambda^2.
$$
We see that we need to find the minimum in $\lambda$ of
\begin{multline*}
\Phi_2\Big((d-1)\lambda^2+\frac{(d-1)^2}{n-d} \lambda^2 \Big)+ z\Big( \Phi_2\rho  + \frac{\Phi_{11}}{2}\Big)- \frac{ (\Phi_{12}z + \frac{(d-1)\Phi_2}{n-d}\lambda)^2}{\Phi_2 \frac{n+1-d}{n-d}+ 2 \Phi_{22}z}\\
=\lambda^2\frac{\Phi_2(d-1)\Big(n\Phi_2+2(n-1)z\Phi_{22}\Big)}{(n-d)(\Phi_2 \frac{n+1-d}{n-d}+ 2 \Phi_{22}z)}
-\lambda\frac{2(d-1)\Phi_2\Phi_{12}z}{(n-d)(\Phi_2 \frac{n+1-d}{n-d}+ 2 \Phi_{22}z)}\\
-\frac{ \Phi_{12}^2z^2}{\Phi_2 \frac{n+1-d}{n-d}+ 2 \Phi_{22}z}
 +z\Big( \Phi_2\rho  + \frac{\Phi_{11}}{2}\Big)
\end{multline*}
Then, once again, to get a minimizer in $\lambda$, we  also
need
$$(n-d)(n\Phi_2+ 2(n-1)\Phi_{22} z )>0,
$$
which is the third assumption in~\eqref{eq-condition-gradient}. 
Then, the minimizer is
\begin{multline*}
\frac{L(\Phi(F, \Gamma(F)))}{2}\geq-\frac{(d-1)\Phi_2\Phi_{12}^2z^2}{(n-d)(n\Phi_2+2(n-1)z\Phi_{22})(\Phi_2 \frac{n+1-d}{n-d}+ 2 \Phi_{22}z)}\\
-\frac{ \Phi_{12}^2z^2}{\Phi_2 \frac{n+1-d}{n-d}+ 2 \Phi_{22}z}+z\Big( \Phi_2\rho  + \frac{\Phi_{11}}{2}\Big)=\frac{-z^2(n-1)\Phi_{12}^2}{n\Phi_2+ 2z(n-1) \Phi_{22}}+ z\Big( \Phi_2\rho  + \frac{\Phi_{11}}{2}\Big)\geq0,
\end{multline*}
which is the fourth assumption in~\eqref{eq-condition-gradient} since $z\geq0$. 

When $\Phi_2=0$ at some point $(y,z)$, we follow the inequality to check the assumptions needed. 
\epf

\bcor[Sub-harmonic functionals under $TCD(0,n)$]
\label{cor-sous-harm-2}
Let assume that $L=\Delta_{\mathfrak{g}}+X$ satisfies a $TCD(0,n)$ with $n<0$. Let $\theta$ be a positive and smooth function on an interval $I\subset \dR$ such that 
\beq
\label{eq-n-admissible}
2\frac{n-1}{n}(\theta')^2\leq \theta\theta''.
\eeq
Then for any  function $F:M\mapsto I$ such that 
$LF=0$, 
\beq
\label{eq-harm-ad}
L(\theta(F)\Gamma(F))\geq0.
\eeq
In particular, for any $\beta\in[n/(2-n),0]$ and nonnegative harmonic function $F$,  
$$
L(F^\beta\Gamma(F))\geq0.
$$
\ecor

When $n$ goes to $-\infty$, the condition~\eqref{eq-n-admissible} on $\theta$ degenerates into $2\theta'\leq \theta\theta''$. This means that the function $\Phi$ such that $\Phi''=\theta$ is admissible : $\Phi$ is convex and $1/\Phi''$ is concave. 

Moreover, when  $n\in(-\infty,0)$, the extremal case is given when $\theta(x)=x^{n/(2-n)}$, where~\eqref{eq-n-admissible} is an equality. 

Although we shall mainly use  Corollary~\ref{cor-sous-harm-2} instead the general result Theorem~\ref{grad.bound.gal.Phi}, the proof in this last case is not really simpler than the general one.

\section{The operators $Q_t^{(m)}$}
\label{sec-Qtm}

\bdefi[The operators $Q_t^{(m)}$]
For any $t\geq0$,  $m>0$, $x\in\bbR^d$, and any bounded  function $f:\bbR^d\mapsto\bbR$, let
\beq
\label{eq-def-q}
Q_t^{(m)}(f)(x)=\frac{1}{c(m,d)}\int \frac{f(yt+x)}{(1+|y|^2)^{\frac{m+d}{2}}}dy,
\eeq
 where 
 \beq
 \label{eq-def-c}
 c(m,d)=\int \frac{1}{(1+|y|^2)^{\frac{m+d}{2}}}dy=\frac{\Gamma(\frac{m}{2})}{\Gamma(\frac{m+d}{2})}\pi^{d/2}
 \eeq
  is the normalization constant, such that $Q_t^{(m)}(1)=1$.
In other words, for any $t>0$, 
$$ 
Q_t^{(m)}(f)(x)=\int f(y)q_t(x,dy),
$$
where the kernel $q_t$ is given by 
\beq
\label{eq-def-kernel}
q_t(x,dy)=\frac{1}{c(m,d)}\frac{t^m}{(t^2+|x-y|^2)^{\frac{m+d}{2}}}dy.
\eeq

\edefi 
The probability  measure
$$
\frac{1}{c(m,d)}(1+|y|^2)^{-\frac{m+d}{2}} dy
$$
in $\dR^d$, with $m>0$,  is known in statistics as the multivariate $t$-distribution with  $m$ degrees of freedom. There is a huge  literature on the subject, see for instance~\cite{kotz}. 

Through an integration by parts (or from the computation of normalization constants) we may notice the useful formula,    
\begin{equation}
\label{eq-ipp-c}
\frac{c(m,d)}{c(m-2,d)}=\frac{m-2}{m-2+d},
\end{equation}
for any $m>2$.

\bprop[Harmonicity of $Q_t^{(m)}f$]
\label{prop-harm}
Let $m>0$. For any smooth and compactly supported function $f$, the map $(0,+\infty)\times\bbR^d\ni(t,x)\mapsto Q_t^{(m)}(f)(x)$ is solution of the elliptic equation $\Delta^{(m)}Q_t^{(m)}(f)=0$ where 
\beq
\label{eq-harm-qt}
\Delta^{(m)}=\Delta_{\bbR^d}+\partial^2_{tt}+\frac{1-m}{t}\partial_t.
\eeq
\eprop
\bpf
This can be proved directly on the expression of the kernel~\eqref{eq-def-kernel}. 
\epf

From Proposition~\ref{prop-harm} (or from a direct computation) the operators $(Q_t^{(m)})_{t\geq 0}$ admit an another formulation.  Let $(Z_s)_{s\geq0}=(X_s,Y_s)_{s\geq0}\in\bbR^{d}\times\bbR$ be a diffusion process with generator $\Delta^{(m)}$ starting from $(x,t)\in\dR^d\times (0,\infty)$. 
That is $(X_s)_{s\geq0}$ is a Brownian  motion  up to a factor $\sqrt{2}$ and $(Y_s)_{s\geq0}$ is independent of $(X_s)_{s\geq0}$ with generator $\partial^2_t+\frac{1-m}{t}\partial_t$.  

The process $(Y_s)_{s\geq0}$ is a Bessel process with parameter $m$ starting from $t>0$. When $m>0$, the stopping time  $S=\inf\{s>0,\,\,Y_s=0\}$ is finite a.s., and 
since $Q_t^{(m)}(f)$ is a solution of the elliptic equation~\eqref{eq-def-q},  from It\^o's Lemma,  or Dynkin's formula~\cite[Sec 7.3]{oksendal}	
\beq
\label{eq-dynkin}
Q_t^{(m)}(f)(x)=E_{x,t}(f(X_S)),
\eeq
for every smooth and bounded function $f$. In other words since, $S$ and $(X_s)_{s\geq0}$ are independent, if $\sigma_m(s,t)ds$ denotes the law of the hitting time $S$ starting from $t>0$, 
\begin{equation}
\label{eq-101}
Q_t^{(m)}(f)(x)=\int_0^\infty P_sf(x)\sigma_m(s,t)ds, 
\end{equation}
where $(P_s)_{s\geq0}$ is the heat semi-group in $\bbR^d$ associated to the Euclidean Laplacian $\Delta$.  From the definition of the operators $Q_t^{(m)}$ and the definition of the the equation $(P_s)_{s\geq0}$, 
$$
P_sf(x)=\int f(y)\exp\Big(-\frac{|x-y|^2}{4s}\Big)\frac{dy}{(4\pi s)^{d/2}},
$$
and this gives a way to compute the law $\sigma_m$ of the hitting time. One obtains,  for any $t>0$, the probability measure
\beq
\label{eq-def-hitting}
\sigma_m(s,t)ds=\frac{1}{2^m\Gamma(m/2)}\frac{t^{m}\exp(-t^2/4s)}{s^{m/2+1}}ds.
\eeq
This formula can also  be found for instance in~\cite{borodin}.  

The map $(t,s)\mapsto \sigma_m(t,s)$ satisfies 
$$
\partial_s=\partial_{tt}^2+\frac{1-m}{t}\partial_t, 
$$ 
with Dirichlet boundary conditions,  $\sigma_m(t,0)=0$ for every $t>0$,  and $\sigma_m(0,s)ds=\delta_0(s)$. 
We recognize in the RHS the generator of the Bessel process. This phenomena is universal and has nothing to do in particular with Bessel processes, see for instance~\cite[Sec 7.2]{pavliotis}.

From this observation and an integration by parts, one recovers that the map $(t,x)\mapsto Q_t^{(m)}f(x)$ satisfies the elliptic equation $\Delta^{(m)}Q_t^{(m)}(f)=0$.

\blem
\label{lem-energie}
Let $m>0$ and  $p\in(0,m/2)$,  then for  any smooth and bounded function $g$, 
$$
E_{x,t}(S^pg(X_S))=t^{2p}\frac{\Gamma(\frac{m}{2}-p)}{4^p\Gamma(\frac{m}{2})}Q_t^{(m-2p)}(g)(x).
$$
In particular, when $p=1$ and $m>2$, 
\beq
\label{eq-energie}
E_{x,t}(Sg(X_S))=\frac{t^{2}}{2(m-2)}Q_t^{(m-2)}(g)(x).
\eeq
\elem

\bpf
The formula comes from a direct computation since we have
\begin{multline*}
E_{x,t}(S^pg(X_S))=\frac{1}{2^m\Gamma(m/2)(4\pi)^{d/2}}\int g(y)\int_0^\infty s^{p-\frac{m+d}{2}-1}{\exp\Big(-\frac{|x-y|^2}{4s}-\frac{t^2}{4s}\Big)}ds\,dy\\
=t^{2p}\frac{\Gamma(\frac{m}{2}-p)}{2^m\Gamma(m/2)(4\pi)^{d/2}}\int P_sg(y)\sigma_{m-2p}(s)ds.
\end{multline*}
\epf

\blem[$QM$ condition for $\Delta^{(m)}$]
\label{lem-tcd-qm}
For any $m\neq 1$, the operator $\Delta^{(m)}$ on $\bbR^{d}\times(0,\infty)$ is a quasi-model $QM(0,d-m+2)$.
\elem

\bpf
For that model, 
$$
\Delta^{(m)}=\Delta_{\bbR^{d+1}}+\frac{1-m}{t}\partial_t,
$$
that is $X=\frac{1-m}{t}\partial_t$. Then $\Ric(L)=\frac{1-m}{t^2}(\partial_t )^2$ and $X\otimes X=\frac{(1-m)^2}{t^2}(\partial_t)^2$ and the result follows from the definition~\eqref{eq-quasi-models}. 
\epf

\section{Beckner and $\Phi$-entropy  inequalities for Cauchy type distribution}
\label{sec-be-cauchy}

\subsection{Beckner inequalities for Cauchy type distribution}
\label{sec-be-cau}
Let $f:\bbR\mapsto [0,\infty)$ be a smooth and compactly supported function,  and let $F(t,x)=Q_t^{(m)}(f)(x)$ with $m>0$.

Let $\beta\in\bbR$.  From proposition~\eqref{prop-harm}, $F$ is $\Delta^{(m)}$-harmonic, and then 
$$
\Delta^{(m)}(F^{\beta+2})=(\beta+2)(\beta+1)F^\beta\Gamma^{\Delta^{(m)}}(F).
$$
and, as in Section~\ref{sec-Qtm}, we can write,  
$$
E_{x,t}(F(Z_S)^{\beta+2})=E_{x,t}(F(Z_0)^{\beta+2})+(\beta+2)(\beta+1)\int_0^\infty E_{x,t}(F^\beta(Z_s)\Gamma^{\Delta^{(m)}}(F)(Z_s)\I_{s\leq S})ds, 
$$
where $(Z_s)_{s\geq0}$ is the Markov process in $\bbR^{d+1}$ with generator $\Delta^{(m)}$.  Since $F(Z_S)=f(X_S)$, from~\eqref{eq-dynkin} we have the general equation,
\beq
\label{eq-ito}
Q_t^{(m)}(f^{\beta+2})=Q_t^{(m)}(f)^{\beta+2}+(\beta+2)(\beta+1)\int_0^\infty E_{x,t}(F^\beta(Z_s)\Gamma^{\Delta^{(m)}}(F)(Z_s)\I_{s\leq S})ds. 
\eeq
Again, from Dynkin's formula applied to $F^\beta\Gamma^{\Delta^{(m)}}(F)$, we have 
\begin{multline*}
E_{x,t}(F^\beta(Z_s)\Gamma^{\Delta^{(m)}}(F)(Z_s)\I_{s\leq S})=\\
E_{x,t}(F^\beta(Z_S)\Gamma^{\Delta^{(m)}}(F)(Z_S)\I_{s\leq S})- E_{x,t}\Big(\int_s^S\Delta^{(m)}(F^{\beta}\Gamma^{\Delta^{(m)}}(F))(Z_u)du\I_{s\leq S}\Big). 
\end{multline*}
In other words, 
\begin{multline}
\label{eq-ito-2}
Q_t^{(m)}(f^{\beta+2})-Q_t^{(m)}(f)^{\beta+2}=(\beta+2)(\beta+1) E_{x,t}\big(SF^\beta(Z_S)\Gamma^{\Delta^{(m)}}(F)(Z_S)\big)\\
-(\beta+2)(\beta+1) \int_0^\infty E_{x,t}\Big(\int_s^S\Delta^{(m)}(F^{\beta}\Gamma^{\Delta^{(m)}}(F))(Z_u)du\I_{s\leq S}\Big)ds. 
\end{multline}

Let now assume that $m\geq d+2$,  and let denote  $n=d-m+2$ ($n\leq0$). Let also assume that  $\beta\in[\frac{n}{2-n},0]$ (that is $\beta\in[-1+\frac{2}{m-d},0]$, $\frac{2}{m-d}>0$). Since $n\leq0$,  $\Delta^{(m)}$ satisfies $TCD(0,n)$ (Lemma~\ref{lem-tcd-qm}) and from Corollary~\ref{cor-sous-harm-2}, 
$$
\Delta^{(m)}(F^{\beta}\Gamma^{\Delta^{(m)}}(F))\geq 0.
$$
The equation~\eqref{eq-ito-2} becomes 
$$
Q_t^{(m)}(f^{\beta+2})- Q_t^{(m)}(f)^{\beta+2}\leq(\beta+2)(\beta+1)\int_0^\infty E_{x,t}(F^\beta(Z_S)\Gamma^{\Delta^{(m)}}(F)(Z_S)\I_{s\leq S})ds,
$$
which is 
$$
Q_t^{(m)}(f^{\beta+2})- Q_t^{(m)}(f)^{\beta+2}\leq(\beta+2)(\beta+1)E_{x,t}(SF^\beta(Z_S)\Gamma^{\Delta^{(m)}}(F)(Z_S)). 
$$
Starting form $Z_0=(x,t)$,  
$$
F^\beta(Z_S)\Gamma^{\Delta^{(m)}}(F)(Z_S)=f^\beta(X_S)\Gamma^{\Delta^{(m)}}(F)(X_S,0)=f^\beta(X_S)\big[\Gamma(f)(X_S)+(\partial_t F)^2(X_S,0)\big],
$$
since $\Gamma^{\Delta^{(m)}}(F)=\Gamma(F)+(\partial_t F)^2$ where $\Gamma(f)=|\nabla f|^2$, ($f:\bbR^d\mapsto \bbR$) is the carr\'e du champ operator associated to Euclidean 
Laplacian in $\bbR^d$. From Proposition~\ref{prop-harm}, since $m>1$, $\partial_t F(x,0)=0$ then  
$$
F^\beta(Z_S)\Gamma(F)(Z_S)=f^\beta(X_S)\Gamma(f)(X_S).
$$
We can now apply Lemma~\ref{lem-energie}, to get
\beq
\label{eq-beckner-qt1}
Q_t^{(m)}(f^{\beta+2})-Q_t^{(m)}(f)^{\beta+2}\leq (\beta+2)(\beta+1)\frac{t^2}{2(m-2)}Q_t^{(m-2)}(f^\beta\Gamma(f)).
\eeq
Changing $f^{\beta+2}$ by $f^2$ and $p=\beta+2$ we have obtained
\bthm[Beckner inequalities for $Q_t^{(m)}$]
\label{thm-beckner-qt}
For any $m\geq d+2$ and $p\in[1+\frac{2}{m-d},2]$, and for any smooth function $f\geq0$,
\beq
\label{eq-beckner-qt}
\frac{p}{p-1}\Big(Q_t^{(m)}(f^2)-Q_t^{(m)}(f^{2/p})^{p}\Big)\leq \frac{2t^2}{m-2}Q_t^{(m-2)}(\Gamma(f)).
\eeq
When $p=2$ ($m\geq d+2$), we obtain the Poincar\'e inequality for every smooth function $f$, 
\beq
\label{eq-poincare-qt}
Q_t^{(m)}(f^2)-Q_t^{(m)}(f)^{2}\leq \frac{t^2}{(m-2)}Q_t^{(m-2)}(\Gamma(f)).
\eeq
\ethm
For $t=1$ and $x=0$ in the previous inequality, we have obtained the Beckner inequality for the Cauchy distribution 
$\frac{1}{c(m,d)}\frac{1}{(1+|y^2|)^{\frac{m+d}{2}}}$.

So, for any $m\geq d+2$ and $p\in[1+\frac{2}{m-d},2]$ one has, 
\begin{multline}
\label{eq-beckner-qt2}
\frac{p}{p-1}\left[\frac{1}{c(m,d)}\int \frac{f^2}{(1+|y|^2)^{\frac{m+d}{2}}}dy-\left(\frac{1}{c(m,d)}\int \frac{f^{2/p}}{(1+|y|^2)^{\frac{m+d}{2}}}dy\right)^{p}\right]\\
\leq \frac{2}{m-2}\frac{1}{c(m-2,d)}\int \frac{\Gamma(f)}{(1+|y|^2)^{\frac{m+d-2}{2}}}dy.
\end{multline}
The inequality can be written as a weighted Poincar\'e-type inequality for a generalized Cauchy distribution. Let define  for $b> \frac{d}{2}$, the probability measure  
$$
d\nu_{b}(y)=\frac{1}{c(2b-d,d)}\frac{1}{(1+|y|^2)^{b}}dy
$$
and replacing  $b=(m+d)/2$ in the previous inequality, we have obtained,  using~\eqref{eq-ipp-c}, 
\bthm[Beckner inequalities for Cauchy distribution]
\label{thm-beckner-cauchy}
For any $b\geq d+1$ and $p\in[1+\frac{1}{b-d},2]$, 
\beq
\label{eq-beckner-cauchy}
\frac{p}{p-1}\left[\int {f^2}d\nu_b-\left(\int {f^{2/p}}d\nu_b\right)^{p}\right]\leq
 \frac{1}{(b-1)}\int \Gamma(f)(1+|y|^2)d\nu_b,
\eeq
for any smooth and nonnegative function $f$.  

Moreover, the inequality is optimal, that is,  fixing  $p\in[1+\frac{1}{b-d},2]$,  the constant $\frac{1}{(b-1)}$ is the best possible in the inequality~\eqref{eq-beckner-cauchy}.
\ethm
The case $p=2$ can be extended to any smooth function $f$ (not only nonnegative function). Then we obtain the following inequality, proved in~\cite{scheffer-phd,blanchet2007,bonforte2010,nguyen-14}.

\bcor[Poincar\'e inequality for Cauchy distribution]
\label{cor-poincare-cauchy}
For any $b\geq d+1$, the measure $\nu_b$ satisfies a Poincar\'e type inequality, 
\beq
\label{eq-poincare-cauchy}
\int {f^2}d\nu_b-\left(\int fd\nu_b\right)^{2}\leq
\frac{1}{2(b-1)}\int \Gamma(f)(1+|y|^2)d\nu_b,
\eeq
for any smooth function $f$. Moreover, functions $y\mapsto y_i$ with $1\leq i\leq d$ are extremal functions. 
\ecor
\bpf
We only have to prove optimality. Apply the Poincar\'e inequality~\eqref{eq-poincare-cauchy} to the function $y\mapsto y_i$ for some $1\leq i\leq d$,
$$
\nu_b(y_i^2)\leq \frac{1}{2(b-1)}\nu_b(1+|y|^2)=\frac{1}{2(b-1)}(1+\nu_b(|y|^2)),
$$  
or equivalently 
$$
\frac{1}{d}\nu_b(|y|^2)\leq \frac{1}{2(b-1)}(1+\nu_b(|y|^2)).
$$
But this last one is an equality since, from~\eqref{eq-ipp-c},  $\nu_b(|y|^2)=\frac{d}{2b-2-d}$.
\epf
\brmqs~
\begin{itemize}
\item Inequality~\eqref{eq-beckner-cauchy} is equivalent to~\eqref{eq-beckner-qt};  indeed we only have to replace  the map  $y\mapsto f(y)$ by   $y\mapsto f(ty+x)$.
\item We cannot reach the logarithmic Sobolev inequality, since in our computation $p\geq1+1/b$ and this inequality would require  a limit $p$ goes to $1$.  
\item A Taylor expansion with $\varepsilon\rightarrow 0$ in inequality~\eqref{eq-beckner-cauchy} with $f=1+\varepsilon g$ implies back the optimal  Poincar\'e inequality~\eqref{eq-poincare-cauchy}, since 
$$
\frac{p}{p-1}\left[\int {f^2}d\nu_b-\left(\int {f^{2/p}}d\nu_b\right)^{p}\right]=2\varepsilon^2\left[\int g^2d\nu_b-\left(\int gd\nu_b\right)^2\right]+o(\varepsilon^2).
$$
We do not know if there are extremal functions apart of constant functions. 
\item Corollary~\ref{cor-poincare-cauchy}  has also been proved by Nguyen~\cite[Cor. 14]{nguyen-14}, using  the Brunn-Minkowski theory (see also the approach in~\cite{bobkov-ledoux, abj}).

In the one dimensional case, the exact value  has been computed by Bonnefont, Joulin and Ma in~\cite[Thm~3.1]{joulin}.  Their result is more general, since they compute the optimal constant for every $b>1/2$. The constant of~\eqref{eq-poincare-cauchy} with $d=1$ is  
$$
\left\{
\begin{array}{l}
\disp\frac{1}{2(b-1)}\quad {\rm if}\quad b\geq3/2;\\
\disp\frac{4}{(2b-1)^2}\quad {\rm if}\quad 1/2< b\leq3/2.
\end{array}
\right. 
$$
It is interesting to notice that there are two regimes, depending on the range of the parameter $b$.  Actually,  we are not able with our method to reach this range,  even in dimension~1. 
\item Applying inequality~\eqref{eq-beckner-cauchy} to the function $x\mapsto f(\sqrt{2b}x)$, then, when $b\rightarrow +\infty$, the inequality becomes the optimal Beckner inequality for the Gaussian measure $\gamma$ in $\dR^d$, 
$$
\frac{p}{p-1}\left[\int {f^2}d\gamma-\left(\int {f^{2/p}}d\gamma\right)^{p}\right]\leq
 2\int \Gamma(f)d\gamma,
$$
for any smooth and positive function $f$. This inequality is proved by W. Beckner~\cite{beckner1}. In that case the inequality holds for any $p\in(1,2]$, and the limit case (when  $p\rightarrow 1$) leads to the optimal Logarithmic Sobolev inequality for the Gaussian measure.  
\item  At  the same time, Nguyen proves  in~\cite{18-nguyen} the Beckner inequality~\eqref{eq-beckner-cauchy} for a range of parameter $p$ strictly contained in our interval $[1+1/(b-d),2]$. Even if, with his method, the parameter $p$ can not reach the full interval, Nguyen is able to extend the result to a general class of probability measures with a convexity assumption.  The method is actually improved in a more general context  in~\cite{18dgz}.
\end{itemize}
\ermqs

\subsection{Tightness of the  inequalities}
\label{sec-tightness}

We are interesting in the tightness of inequality~\eqref{eq-beckner-qt},  in the sense of the curvature-dimension condition. 
In all this section we assume that $d\geq 2$.

\blem[Taylor expansion of $(Q_t^{(m)}$]
\label{lem-taylor-ex}
For any $m>4$, we have 
\beq
\label{eq-taylor-ex}
Q^{(m)}_t=\Id+\frac{\Delta}{2(m-2)}t^2+\frac{\Delta^2}{8(m-2)(m-4)}t^4+o(t^4).
\eeq
 The formula has to be understood as follow, for any smooth and compactly supported function~$f$, the rest function $o(t^4)$ is uniformly bounded in the variable $x$. 
\elem
\bpf
One can check directly on the formula. On the other hand, one can deduce the formula with functional analysis, at least formally.  Since $Q_t^{(m)}$ is $\Delta^{(m)}$-harmonic from Proposition~\ref{prop-harm} then, from the functional analysis point of view, $Q_t^{(m)}$ is a function $G_m(t,\Delta)$. From the definition~\eqref{eq-harm-qt} of $\Delta^{(m)}$, the function $G_m$ satisfies for  $t\geq0$, $x<0$, 
$$
xG_m(t,x)+\frac{1-m}{t}\partial_tG_m(t,x)+\partial_{tt}^2G_m(t,x)=0, 
$$  
with the boundary conditions $G_m(0,x)=1$. Solution of such differential equation is given by, $G_m(t,x)=H_m(-t^2x)$ where $H_m$ satisfies for any $\lambda\geq0$, 
$$
4\lambda H_m''(\lambda)-2(m-2)H'_m(\lambda)-H_m(\lambda)=0.
$$
An asymptotic development of $H_m(\lambda)=\sum_{p\geq0}c_p\lambda^p$ gives $c_0=1$ (since $G_m(0,x)=1$) and 
$$
c_{p+1}=\frac{1}{2(p+1)(2p-m+2)}c_p.
$$
For the forth order Taylor expansion of $Q_t^{(m)}$,  we get~\eqref{eq-taylor-ex}, that is  $c_1=\frac{1}{2(2-m)}$ and $c_2=\frac{1}{8(m-4)(m-2)}$. 
\epf

This Lemma is fundamental. The Taylor expansion in~\eqref{eq-beckner-qt}, or equivalently in~\eqref{eq-poincare-qt}, gives in return the $CD(0,d)$ conditions. We say then that the inequalities are $\Gamma_2$-tight and nothing has been lost in term of curvature-dimension condition~\eqref{eq-def-gamma2}, in other words the Taylor expansion of the inequality implies the $CD(0,d)$ condition. 

To see that, it is enough to compute the Taylor expansion of the inequality~\eqref{eq-beckner-qt},  when $t$ goes to 0. It is important to fix $\beta=\frac{2-m+d}{m-d}=-1+\frac{2}{m-d}$ the extremal bound of the admissible interval of $\beta$.

It is clear that the zero order terms vanish, it is also not difficult to observe that the second order terms also vanish. Therefor we only have to compute the fourth order terms and we obtain the following expression, 
\begin{multline*}
\frac{1}{8(m-4)(m-2)}\Delta^2(f^{\beta+2})-\frac{(\beta+2)}{8(m-4)(m-2)}f^{\beta-1}\Delta^2(f)-\frac{(\beta+2)(\beta+1)}{8(m-2)^2}f^{\beta-1}\Delta(f)^2\leq\\
\frac{\beta(\beta+1)(\beta+2)}{2(m-4)(m-2)}\Gamma(f,\Gamma(f))f^{\beta-1}+\frac{\beta(\beta+1)(\beta+2)}{4(m-4)(m-2)}f^{\beta-1}\Delta (f)\Gamma(f)\\
+\frac{(\beta-1)\beta(\beta+1)(\beta+2)}{4(m-4)(m-2)}f^{\beta-2}\Gamma(f)^2+\frac{(\beta+1)(\beta+2)}{4(m-4)(m-2)}f^{\beta}\Delta \Gamma(f).
\end{multline*}
Expanding the term $\Delta^2(f^{\beta+2})$ (recall that $\beta=\frac{2-m+d}{m-d}$), after a bit of algebra, we  finally get, 
\beq
\label{eq-cd-1}
\Gamma_2(f)\geq \frac{\beta+1}{d(\beta+1)-2\beta}(\Delta f)^2-\beta\frac{\Gamma(f,\Gamma(f))}{f}-\frac{\beta(\beta-1)}{2}\frac{\Gamma(f)^2}{f^2}.
\eeq
First, observe that if $\beta=0$ that is $m=d+2$ then the previous inequality becomes 
$$
\Gamma_2(f)\geq \frac{1}{d}(\Delta f)^2,
$$ 
in others words we recover the $CD(0,d)$ condition of the Laplacian $\Delta$. 

It is more subtle  to notice that, for any $\beta\in(-1,0)$ fixed, inequality~\eqref{eq-cd-1} implies the same condition  $CD(0,d)$.  It is more tricky since, taking $f=1+\varepsilon g$ with $\varepsilon\rightarrow 0$, then~\eqref{eq-cd-1} already implies $\Gamma_2(g)\geq \frac{\beta+1}{d(\beta+1)-2\beta}(\Delta g)^2$ which is  the $CD(0,\frac{d(\beta+1)-2\beta}{\beta+1})$ condition, weaker than the  $CD(0,d)$ condition,  since  $\frac{d(\beta+1)-2\beta}{\beta+1}\geq d$. The way to recover the $CD(0,d)$ condition is to replace $f$ by $\Phi(f)$ in~\eqref{eq-cd-1}, where $\Phi$ is a smooth function. 

The diffusion properties in~\eqref{eq-cd-1} one has, 
\begin{multline*}
\Phi'^2(\Gamma_2(f)-a(\Delta f)^2)+\Phi''\Gamma(f)^2(1-a)-c\frac{\Phi'^4}{\Phi^2}\Gamma(f)^2+\Phi'\Phi''(\Gamma(f,\Gamma(f))-2a\Delta f\Gamma(f))\\
-2b\frac{\Phi'^2\Phi''}{\Phi}\Gamma(f)^2-b\frac{\Phi'^3}{\Phi}\Gamma(f,\Gamma(f))\geq0,
\end{multline*}
where $a=\frac{\beta+1}{d(\beta+1)-2\beta}$, $b=-2\beta$ and $c=-\beta(\beta-1)$ and for simplicity we omit the variable $f$ in the function $\Phi$. 
It implies that the quadratic form with respect to the variables $\Phi'$, $\Phi'^2/\Phi$ and $\Phi''$ is positive, then the determinant of its matrix is positive. We get 
\begin{multline*}
\Gamma_2(f)\geq \frac{1}{d}(\Delta f)^2+\PAR{\frac{-ac-a^2}{c(a-1)-b^2}-\frac{1}{d}}(\Delta f)^2+\frac{-b^2-c-b^2a}{4(c(a-1)-b^2)}\PAR{\frac{\Gamma(f,\Gamma(f))}{\Gamma(f)}}^2\\
+\frac{b^2+ac}{c(a-1)-b^2}\frac{\Gamma(f,\Gamma(f))}{\Gamma(f)}\Delta f.
\end{multline*}
After an unpleasant computation, the inequality is independent of $\beta$, and may be written as follow
$$
\Gamma_2(f)\geq \frac{1}{d}(\Delta f)^2+\frac{d}{d-1}\PAR{\frac{\Gamma(f,\Gamma(f))}{2\Gamma(f)}-\frac{1}{d}\Delta f}^2,
$$
which is a reinforced $CD(0,d)$ condition. 

Indeed, in all these computations, we could have replaced the Laplace operator $\Delta$ in $\bbR^d$ by any operator satisfying the $CD(0,d)$ condition. Then we could have  constructed the associated $Q_t^{(m)}$ operator through~\eqref{eq-101}, and obtain also a similar family of Beckner inequalities for this $Q_t^{(m)}$, still equivalent to the starting  $CD(0,d)$ condition. The main difference in the Euclidean case is that then the family of inequalities for $Q_t^{(m)}(x)$  reduces one for $Q_1^{(m)}(0)$, through dilations and translations.  
\subsection{$\Phi$-entropy inequalities for Cauchy type distribution}
\label{sec-phi-cau}

We can extend Beckner's inequalities to a general $\Phi$-entropy inequality. For any convex function $\Phi$ on an interval $I\subset \bbR$, let  us define the $\Phi$-entropy for  measure $\mu$ (or a kernel) and a function $f$ on $I$,  
$$
{\rm Ent}_{\mu}^\Phi(f)=\int \Phi(f)d\mu-\Phi\Big(\int fd\mu\Big). 
$$
Section~\ref{sec-be-cau} can be easily generalized to $\Phi$-entropy. 
\bdefi[$n$-Admissible functions]
\label{def-n-nadmissible}
Let $\Phi:I\mapsto \dR$ be a  smooth function on an interval $I\subset \dR$ and let $n<0$. We say that $\Phi$ is $n$-admissible if 
$\Phi''>0$ on $I$ and $(\Phi'')^{\frac{2-n}{n}}$ is concave. 

In other words, $\Phi$ is $n$-admissible if and only if  $\Phi''>0$  and moreover
$$
2\frac{n-1}{n}(\Phi^{(3)})^2\leq \Phi''\Phi^{(4)},
$$ 
that is $\theta=\Phi''$ satisfies condition~\eqref{eq-n-admissible}. 
\edefi

When $n$ goes to $-\infty$, an $(-\infty)$-admissible is just a function $\Phi$ such that,  $\Phi''>0$ on $I$ and $(1/\Phi'')$ is concave. Many $\Phi$-entropy inequalities  have been proved  for such function $\phi$ under the curvature-dimension condition $CD(\rho,\infty)$ (see~\cite{chafai04,bolley-gentil}). In our computations, the case $n=-\infty$ is then similar to the case $n=+\infty$. 

We can state the following result.

\bthm[$\Phi$-entropy inequality for $Q_t^{(m)}$]
For any function $\Phi$ $(d-m+2)$-admissible on an interval $I$, the following holds, for any $m\geq d+2$,
\beq
\label{eq-phi-ent}
{\rm Ent}_{Q_t^{(m)}}^\Phi(f)\leq \frac{t^2}{2(m-2)}Q_t^{(m-2)}(\Phi''(f)\Gamma(f)),
\eeq
for any smooth function $f$ on $I$. 

The inequality is optimal in the sense that if inequality~\eqref{eq-phi-ent} holds for a some function $\Phi$, the constant $\frac{t^2}{2(m-2)}$ is the optimal one. \ethm

\bcor[$\Phi$-entropy inequality for Cauchy distribution]
For any function $\Phi$ $(d-m+2)$-admissible on an interval $I$, the following holds, for $b\geq d+1$
$$
{\rm Ent}_{\nu_b}^\Phi(f)\leq \frac{1}{4(b-1)}\int \Phi''(f)\Gamma(f)(1+|x|^2)d\nu_b,
$$
for any smooth function $f$ on $I$.  The inequality is optimal. 
\ecor

The proof is the same as  the one of Theorem~\ref{thm-beckner-qt}, using the sub-harmonicity inequality~\eqref{eq-harm-ad}, so we skip it, the optimality being proved in a same  way. 

\section{A Beckner-type inequality on the sphere}
\label{sec-be-sphere}
In this section we are looking for the link between the family of  Beckner inequalities~\eqref{thm-beckner-cauchy} for $\nu_b$ (or $Q_t^{(m)}$) and a family of  Beckner inequalities on the sphere. In all this section we assume that $d\geq 2$.

The unit sphere $\bbS^d\subset\bbR^{d+1}$ can be seen in $\bbR^d$ through the stereographic projection, with the carr\'e du champ operator on the sphere 
$$
\Gamma_{\bbS}(f)=\frac{\rho^4}{4}\Gamma(f),
$$ 
where $\rho(x)=\sqrt{1+|x|^2}$ and the spherical measure 
$$
\mu_{\bbS}(dx)=\frac{1}{c(d,d)}\frac{1}{\rho^{2d}}dx,
$$ 
where the normalization constant $c(d,d)$ has been defined in~\eqref{eq-def-c} (see~\cite[Sec~2.2]{bgl-book}). The Laplace-Beltrami operator takes the form in the stereographic projection representation 
$$
\Delta_\bbS=\frac{(1+|x|^2)^2}{4}\Delta-\frac{d-2}{2}(1+|x|^2)\sum_{i=1}^dx_i\partial_i,
$$
where $\Delta$ is the classical Laplacian in $\bbR^d$. The operator is symmetric in $L^2(\mu_{\bbS})$.  

It is interesting to recall that the map $\bbS^d\ni x\mapsto x_{d+1}$ is a eigenvector associated to the eigenvalue $-d$. This function, from the stereographic projection, takes the form 
\beq
\label{eq-def-u}
u(x)=\frac{1-|x|^2}{1+|x|^2}, \,\, x\in\bbR^d.
\eeq
Then, $u$ satisfies $\Delta_\bbS u=-du$ and $\Gamma_{\bbS^d}(u)=1-u^2$ since,
$$
\Gamma_{\bbS^d}(u)=\frac{(1+|x|^2)^2}{4}\Gamma\PAR{\frac{1-|x|^2}{1+|x|^2}}=\frac{1}{(1+|x|^2)^2}\Gamma(|x|^2)=\frac{4|x|^2}{(1+|x|^2)^2}=1-u^2.
$$
Indeed this computation is much easer (and indeed completely elementary)  using  the standard representation in the sphere, see~\cite[Sec~2.2]{bgl-book}.

\medskip

We can now transform the Beckner inequality for $Q_t^{(m)}$ into a new Beckner inequality on the sphere.  
With these new notations, inequality~\eqref{eq-beckner-qt1} for $t=1$ and $x=0$  becomes for smooth and nonnegative function $f$,
\begin{multline}
\label{eq-beckner-sphere1}
\frac{c(d,d)}{c(m,d)}\int f^{\beta+2}\rho^{d-m}d\mu_\bbS-\left(\frac{c(d,d)}{c(m,d)}\int f\rho^{d-m}d\mu_\bbS\right)^{\beta+2}\\
\leq\frac{2(\beta+1)(\beta+2)}{m-2}\frac{c(d,d)}{c(m-2,d)}\int f^\beta\Gamma_\bbS(f)\rho^{d-m-2}d\mu_\bbS.
\end{multline}
We concentrate to the case $\beta=\frac{d+2-m}{m-d}$, the limit of  the admissible interval. Let $f=\rho^{m-d}g$ and we compute the various terms. First, we have form diffusion properties of $\Gamma_\bbS$,  
\begin{multline*}
f^\beta\Gamma_\bbS(f)\rho^{d-m-2}=g^\beta\rho^{2(d-m)}\Gamma_\bbS(\rho^{m-d}g)=\\
g^\beta\Gamma_\bbS(g)+2(m-d)g^{\beta+1}\Gamma_\bbS(\log\rho,g)+(m-d)^2g^{\beta+2}\Gamma_\bbS(\log\rho),
\end{multline*}
so  that
\begin{multline*}
\int f^\beta\Gamma_\bbS(f)\rho^{d-m-2}d\mu_\bbS=\\
\int g^\beta\Gamma_\bbS(g)d\mu_\bbS+(m-d)\int \frac{2}{\beta+2}\Gamma_\bbS(\log\rho,g^{\beta+2})d\mu_\bbS+(m-d)^2\int g^{\beta+2}\Gamma_\bbS(\log\rho)d\mu_\bbS.
\end{multline*}
With an integration by parts, we get
$$
\int \Gamma_\bbS(\log\rho,g^{\beta+2})d\mu_\bbS=-\int g^{\beta+2}\Delta_\bbS \log\rho d\mu_\bbS,
$$
and then
$$
\int f^\beta\Gamma_\bbS(f)\rho^{d-m-2}d\mu_\bbS=\int g^\beta\Gamma_\bbS(g)d\mu_\bbS+\int g^{\beta+2}Kd\mu_\bbS,
$$
where 
$$
K=(m-d)^2\Big(2\frac{\Delta_\bbS \log \rho}{d-m-2}+\Gamma_\bbS(\log\rho)\Big),
$$
since $\beta+2=(d-m-2)/(d-m)$. Secondly we have 
$$
\int f^{\beta+2}\rho^{d-m}d\mu_\bbS=\int g^{\beta+2}\rho^{2}d\mu_\bbS
$$
and 
$$
\int f\rho^{d-m}d\mu_\bbS=\int gd\mu_\bbS.
$$
Then, inequality~\eqref{eq-beckner-sphere1} can be written as 
\begin{multline*}
\frac{c(d,d)}{c(m,d)}\int g^{\beta+2}\rho^{2}d\mu_\bbS-\left(\frac{c(d,d)}{c(m,d)}\int gd\mu_\bbS\right)^{\beta+2}\\
\leq\frac{2(\beta+1)(\beta+2)}{m-2}\frac{c(d,d)}{c(m-2,d)}\Big(\int g^\beta\Gamma_\bbS(g)d\mu_\bbS+\int g^{\beta+2}Kd\mu_\bbS\Big),
\end{multline*}
or 
$$
\int g^{\beta+2}Rd\mu_\bbS-\left(\frac{c(d,d)}{c(m,d)}\int gd\mu_\bbS\right)^{\beta+2}\leq\frac{2(\beta+1)(\beta+2)}{m-2}\frac{c(d,d)}{c(m-2,d)}\int g^\beta\Gamma_\bbS(g)d\mu_\bbS,
$$
with
$$
R=\frac{c(d,d)}{c(m,d)}\rho^2-\frac{2(\beta+1)(\beta+2)}{m-2}\frac{c(d,d)}{c(m-2,d)}K.
$$
From~\eqref{eq-ipp-c}, the definition of $\beta$ and the fact that $\rho^2=2/(1+u)$ where $u$ has been defined in~\eqref{eq-def-u}, we can write 
$$
R=\frac{c(d,d)}{c(m,d)}\PAR{\frac{2}{1+u}-2\frac{m-d+2}{(m-d)^2}\frac{K}{m-2+d}}, 
$$
and the inequality 
$$
\int g^{\beta+2}Rd\mu_\bbS-\left(\frac{c(d,d)}{c(m,d)}\int gd\mu_\bbS\right)^{\beta+2}\leq\frac{4(m-d+2)}{(m-d)^2(m-2+d)}\frac{c(d,d)}{c(m,d)}\int g^\beta\Gamma_\bbS(g)d\mu_\bbS.
$$

We need now to compute $\Delta_\bbS\log \rho$ and  $\Gamma_\bbS(\log\rho)$. We have 
\begin{multline*}
\Delta_\bbS(\log \rho)=\frac{1}{2}\Delta_\bbS\log \rho^2=-\frac{1}{2}\Delta_\bbS \log (1+u)=\frac{d}{2(1+u)}u+\frac{1}{2(1+u)^2}\Gamma_{\bbS^d}(u)\\
=\frac{d}{2(1+u)}u+\frac{1}{2(1+u)^2}(1-u^2). 
\end{multline*}
Then, finally 
$$
\Delta_\bbS(\log \rho)=\frac{1+u(d-1)}{2(1+u)},
$$
and
$$
\Gamma_\bbS(\log\rho)=\frac{1-u}{4(1+u)}.
$$
If we plug those values on the definition of $K$, then  a miracle occurs. The function $R$ is constant,
$$
R=\frac{c(d,d)}{c(m,d)}\frac{3d+m-2}{d+m-2}. 
$$
Setting now $f^2=g^{\beta+2}=g^\frac{d-m-2}{d-m}$, we obtain 
\bthm[Beckner-type inequalities on the sphere]
\label{thm-beckner-sphere}
For any smooth nonnegative function $f$ on the sphere $\bbS^d$, and $m\geq d+2$, we have 
\beq
\label{eq-beckner-sphere2}
\int f^{2}d\mu_\bbS\leq A\left(\int f^{2/p}d\mu_\bbS\right)^{p}
+\frac{16}{(m+2-d)(3d-2+m)}\int \Gamma_\bbS(f)d\mu_\bbS,
\eeq
where 
$$
p=1+\frac{2}{m-d}\in(1,2]
$$
and
$$
A=\left(\frac{c(d,d)}{c(m,d)}\right)^{\frac{2}{m-d}}\frac{m+d-2}{m+3d-2}.
$$
\ethm

\brmqs
~

\begin{itemize}
\item The function $R$ is constant only for the parameter $p=1+\frac{2}{m-d}$.
\item  When $m>d+2$, $A>1$. In that case inequality~\eqref{eq-beckner-sphere2} is not tight, in the sense that when $g=1$, we do not have an equality. But the inequality is still optimal since it is saturated for $f=\rho^{\frac{d-m-2}{2}}$.  

\item When $m=d+2$, then $A=1$, $p=2$ and inequality~\eqref{eq-beckner-sphere2} is nothing else than the optimal Poincar\'e inequality on the sphere,
$$
\int f^{2}d\mu_\bbS-\left(\int f^{2}d\mu_\bbS\right)^{2}
\leq\frac{1}{d}\int \Gamma_\bbS(f)d\mu_\bbS.
$$
\item The classical Beckner inequalities hold  for the spherical model (see e.g.~\cite[Rmk 6.8.4]{bgl-book}):  for any $p\in(1,2]$ 
\beq
\label{eq-class-b}
\int f^{2}d\mu_\bbS\leq \left(\int |f|^{2/p}d\mu_\bbS\right)^{p}
+\frac{p-1}{pd}\int \Gamma_\bbS(f)d\mu_\bbS.
\eeq

Inequality~\eqref{eq-class-b} is optimal in the sense that $\frac{p-1}{pd}$ is the best constant. Only constant functions saturate 
inequality~\eqref{eq-class-b} for $p\in(1,2)$.  And when $p$ goes to 1, inequality~\eqref{eq-class-b} provides the Logarithmic Sobolev inequality on the sphere
$$
\int f^{2}\log \frac{f^{2}}{\int f^{2}d\mu_\bbS}d\mu_\bbS\leq \frac{2}{d}\int \Gamma_\bbS(f)d\mu_\bbS,
$$
proved in~\cite{weissler} (see also~\cite[Thm.~5.7.4]{bgl-book}). On the other hand, this last logarithmic Sobolev inequality implies back the full family~\eqref{eq-class-b}, $p\in(1,2]$.

Let us mention  that inequality~\eqref{eq-class-b} has been generalized by in~\cite{dekl} in the following way,  
$$
\Psi\left( \frac{1-\int f^{2/p}d\mu_\bbS}{p-1}\right)\leq \int \Gamma_\bbS(f)d\mu_\bbS, 
$$
for every nonnegative function $f$ such that $\int f^2d\mu_\bbS=1$ and for some explicit function $\Psi$.

\item So~\eqref{eq-beckner-sphere2} appears as a new and optimal inequality on the sphere. It is probably worth to compare it to the del Pino-Dolbeault family of optimal Gagliardo-Nirenberg inequalities, 
$$
\|f\|_{2\frac{a-1}{a-2}}\leq C\,\|\nabla f \|_{2}^\theta\,   \|f\|_{\frac{2a}{a-2}}^{1-\theta},\quad a\geq d,
$$
where $\theta$ is a fixed by scaling properties, see~\cite{pino-dolbeault}.
\item The  Taylor expansion (when $m\rightarrow\infty$) of the constant $A$ is given by 
$$
A=1+d\frac{\log(m)}{m}+\frac{C}{m}+o(1/m),
$$
where $C$ is a constant depending on the dimension $d$. This will be used in the rest of the section. 
\end{itemize}
\ermqs

Even if we are not able from Theorem~\ref{thm-beckner-sphere} to reach directly the Logarithmic Sobolev inequality,  inequality~\eqref{eq-beckner-sphere2} contains enough information to obtain a Sobolev inequality for the spherical model, however with a non optimal constant. 

For this, we are going to obtain from~\eqref{eq-beckner-sphere2} a Nash inequality on the sphere. 

First, it is not difficult to see that there exits two constants $\alpha, \beta>0$ depending only on $d$, such that, for any $m\geq d+2$, 
$$
A\leq \alpha^{\frac{2}{m-d}}m^{\frac{d}{m-d}},
$$
which is just a precise form of the previous Taylor expansion,  
and 
$$
\frac{16}{(m+2-d)(3d-2+m)}\leq \frac{\beta}{m^2}. 
$$
From H\"older's inequality, if $p\in[1,2]$, 
$$
\left(\int |f|^{2/p}d\mu_\bbS\right)^p\leq \left(\int |f|d\mu_\bbS\right)^{p-1}\left(\int f^{2}d\mu_\bbS\right)^{2-p}
$$
with $p=1+2/(m-d)\in[1,2]$. Inequality~\eqref{eq-beckner-sphere2} becomes, for any nonnegative function $f$ such that $\int fd\mu_\bbS=1$, and any $m\geq m+2$,  
$$
\int f^{2}d\mu_\bbS\leq  \alpha^{\frac{2}{m-d}}m^{\frac{d}{m-d}} \left(\int f^{2}d\mu_\bbS\right)^{1-\frac{2}{m-d}}+\frac{\beta}{m^2}\int \Gamma_\bbS(f)d\mu_\bbS. 
$$
Let define for $x\geq0$, the map 
$$
\phi_m(x)=\alpha^{\frac{2}{m-d}}m^{\frac{d}{m-d}}x^{1-\frac{2}{m-d}}+\frac{\beta}{m^2}E,
$$
where $E=\int \Gamma_\bbS(f)d\mu_\bbS$. The function $\phi_m$ as a unique fixed point $x_m>0$, and the previous inequality implies that  $\int f^{2}d\mu_\bbS\leq x_m$.

Let us prove that for any $m\geq d+2$, 
$$
x_m\leq \alpha \,m^{d/2}\left(1+\frac{\beta E}{2\alpha m^{1+d/2}}\right)=a_m. 
$$
To prove such an inequality, it is enough to prove that, for any $m\geq d+2$,  $\phi_m(a_m)\leq a_m$. We have 
\begin{multline*}
\phi_m(a_m)=\alpha \,m^{d/2} \left(1+\frac{\beta E}{2\alpha m^{1+d/2}}\right)^{1-\frac{2}{m-d}}+\frac{\beta}{m^2}E\\
\leq \alpha \,m^{d/2}\left[ \left(1+\frac{\beta E}{2\alpha m^{1+d/2}}\right)^{1-\frac{2}{m}}+\frac{\beta}{m^{2+\frac{d}{2}}}E\right]
\leq \alpha \,m^{d/2}\left(1+\frac{\beta E}{2\alpha m^{1+d/2}}\right)=a_m,
\end{multline*}
where we used the inequality $(1+x)^u\leq 1+ux$, for $x>0$ and $u\in(0,1)$. Modifying the numerical constants (depending only on $d$), we see that for any $m\geq d+2$, 
$$
\int f^{2}d\mu_\bbS\leq C m^{d/2}+\frac{C}{m^{d/2}}\int \Gamma_\bbS(f)d\mu_\bbS.
$$
Optimizing with respect to the parameter $m\geq d+2$, we have obtained, for some other constant $C>0$, 
\begin{multline*}
\int f^{2}d\mu_\bbS\leq C\left(\left(\int |f|d\mu_\bbS\right)^{2}
+\int \Gamma_\bbS(f)d\mu_\bbS\right)^{\frac{d}{d+2}}\left(\int |f|d\mu_\bbS\right)^{\frac{2}{d+2}}\\
\leq C\left(\int f^2d\mu_\bbS
+\int \Gamma_\bbS(f)d\mu_\bbS\right)^{\frac{d}{d+2}}\left(\int |f|d\mu_\bbS\right)^{\frac{2}{d+2}},
\end{multline*}
from Jensen inequality.

This last inequality is a so-called Nash inequality, and is known to be equivalent to the Sobolev inequality (up to the optimal constant) for the spherical model, see~\cite[Prop.~6.2.3]{bgl-book}, 
$$
\left(\int | f|^{\frac{2d}{d-2}}d\mu_\bbS\right)^{\frac{d-2}{d}}\leq C\int f^2d\mu_\bbS +C\int \Gamma_\bbS(f)d\mu_\bbS.
$$
 From the Poincar\'e inequality ($m=d+2$) and~\cite[Prop. 6.2.2]{bgl-book}, this inequality implies a tight Sobolev inequality. We have obtained
\bthm[From Beckner inequalities to Sobolev inequality on the sphere]
\label{thm-sob-sphere}
From the family of Beckner type inequalities~\eqref{eq-beckner-sphere2}, one can prove that there exists a constant $C>0$ depending only on $d$ such that, for any smooth  functions $f$ on the sphere $\bbS^d$, 
\beq
\label{eq-sob-sphere}
\left(\int |f|^{\frac{2d}{d-2}}d\mu_\bbS\right)^{\frac{d-2}{d}}\leq \int f^2d\mu_\bbS +C\int \Gamma_\bbS(f)d\mu_\bbS.
\eeq
\ethm

\bigskip

\noindent 
{\bf Conclusion:} The family of optimal inequalities~\eqref{eq-beckner-sphere2} appears as a new form of optimal inequalities on the sphere,  implying a Sobolev inequality. 

\bigskip

\bigskip

\noindent{\bf Acknowledgements.}   
 This work was partly written while the two first  authors were visiting Institut Mittag-Leffler in Stockholm; it is a pleasure for them to thank this institution for its kind hospitality. We would like also to warmly thank the referee who give interesting comments and pointing out references.
 
 This research was supported  by the French ANR-17-CE40-0030 EFI project.
\medskip


\begin{thebibliography}{{{\O}k}03}

\bibitem[ABJ18]{abj}
M.~{Arnaudon}, M.~{Bonnefont}, and A.~{Joulin}.
\newblock {Intertwinings and generalized Brascamp-Lieb inequalities.}
\newblock {\em {Rev. Mat. Iberoam.}}, 34(3):1021--1054, 2018.

\bibitem[AD05]{arnolddolbeault05}
A.~Arnold and J.~Dolbeault.
\newblock Refined convex {S}obolev inequalities.
\newblock {\em J. Funct. Anal.}, 225(2):337--351, 2005.

\bibitem[Bak94]{bakrystflour}
D.~Bakry.
\newblock L'hypercontractivit\'e et son utilisation en th\'eorie des
  semigroupes.
\newblock In {\em Lectures on probability theory ({S}aint-{F}lour, 1992)},
  Lecture Notes in Math. 1581, pages 1--114. Springer, Berlin, 1994.

\bibitem[BBD{\etalchar{+}}07]{blanchet2007}
A.~Blanchet, M.~Bonforte, J.~Dolbeault, G.~Grillo, and J.-L. V{\'a}zquez.
\newblock Hardy-{P}oincar\'e inequalities and applications to nonlinear
  diffusions.
\newblock {\em C. R. Math. Acad. Sci. Paris}, 344(7):431--436, 2007.

\bibitem[BCLS95]{bcls}
D.~{Bakry}, T.~{Coulhon}, M.~{Ledoux}, and L.~{Saloff-Coste}.
\newblock {Sobolev inequalities in disguise.}
\newblock {\em {Indiana Univ. Math. J.}}, 44(4):1032--1074, 1995.

\bibitem[BDGV10]{bonforte2010}
M.~Bonforte, J.~Dolbeault, G.~Grillo, and J.~L. Vazquez.
\newblock Sharp rates of decay of solutions to the nonlinear fast diffusion
  equation via functional inequalities.
\newblock {\em Proc. Natl. Acad. Sci. USA}, 107(38):16459--16464, 2010.

\bibitem[B{\'E}85]{bakryemery}
D.~Bakry and M.~{\'E}mery.
\newblock Diffusions hypercontractives.
\newblock In {\em S\'eminaire de probabilit\'es, {XIX}, 1983/84}, Lecture Notes
  in Math. 1123, pages 177--206. Springer, Berlin, 1985.

\bibitem[{Bec}89]{beckner1}
W.~{Beckner}.
\newblock {A generalized Poincar\'e inequality for Gaussian measures.}
\newblock {\em {Proc. Am. Math. Soc.}}, 105(2):397--400, 1989.

\bibitem[BG10]{bolley-gentil}
F.~{Bolley} and I.~{Gentil}.
\newblock {Phi-entropy inequalities for diffusion semigroups.}
\newblock {\em {J. Math. Pures Appl. (9)}}, 93(5):449--473, 2010.

\bibitem[BGL14]{bgl-book}
D.~{Bakry}, I.~{Gentil}, and M.~{Ledoux}.
\newblock {\em {Analysis and geometry of Markov diffusion operators.}}
\newblock Cham: Springer, 2014.

\bibitem[BJM16]{joulin}
M.~Bonnefont, A.~Joulin, and Y.~Ma.
\newblock A note on spectral gap and weighted {P}oincar\'e inequalities for
  some one-dimensional diffusions.
\newblock {\em ESAIM Probab. Stat.}, 20:18--29, 2016.

\bibitem[BL09]{bobkov-ledoux}
S.~G. {Bobkov} and M.~{Ledoux}.
\newblock {Weighted Poincar\'e-type inequalities for Cauchy and other convex
  measures.}
\newblock {\em {Ann. Probab.}}, 37(2):403--427, 2009.

\bibitem[BS96]{borodin}
A.N. {Borodin} and P.~{Salminen}.
\newblock {\em {Handbook of Brownian motion - facts and formulae.}}
\newblock Basel: Birkh\"auser, 1996.

\bibitem[Cha04]{chafai04}
D.~Chafa{\"{\i}}.
\newblock Entropies, convexity, and functional inequalities: on
  {$\Phi$}-entropies and {$\Phi$}-{S}obolev inequalities.
\newblock {\em J. Math. Kyoto Univ.}, 44(2):325--363, 2004.

\bibitem[DD02]{pino-dolbeault}
M.~{Del Pino} and J.~{Dolbeault}.
\newblock {Best constants for Gagliardo-Nirenberg inequalities and applications
  to nonlinear diffusions.}
\newblock {\em {J. Math. Pures Appl. (9)}}, 81(9):847--875, 2002.

\bibitem[DEKL14]{dekl}
J.~{Dolbeault}, M.~J. {Esteban}, M.~{Kowalczyk}, and M.~{Loss}.
\newblock {Improved interpolation inequalities on the sphere.}
\newblock {\em {Discrete Contin. Dyn. Syst., Ser. S}}, 7(4):695--724, 2014.

\bibitem[DGS18]{18dgz}
D.~Dupaigne, I.~Gentil, and Zugmeyer S.
\newblock A family of {B}eckner inequalities under different
  curvature-dimension conditions.
\newblock 2018.

\bibitem[DNS08]{dns}
J.~Dolbeault, B.~Nazaret, and G.~Savar\'e.
\newblock {On the Bakry-Emery criterion for linear diffusions and weighted
  porous media equations.}
\newblock {\em Comm. Math. Sci.}, 6(2):477--494, 2008.

\bibitem[IV17]{17-volberg}
P.~{Ivanisvili} and A.~{Volberg}.
\newblock {Improving Beckner's bound via Hermite functions.}
\newblock {\em {Anal. PDE}}, 10(4):929--942, 2017.

\bibitem[IV18]{18Ivanisvili-volberg}
P.~Ivanisvili and A.~Volberg.
\newblock Isoperimetric functional inequalities via the maximum principle: the
  exterior differential systems approach.
\newblock In {\em 50 years with {H}ardy spaces}, volume 261 of {\em Oper.
  Theory Adv. Appl.}, pages 281--305. Birkh\"{a}user/Springer, Cham, 2018.

\bibitem[KN04]{kotz}
S.~{Kotz} and S.~{Nadarajah}.
\newblock {\em {Multivariate $t$ distributions and their applications.}}
\newblock Cambridge University Press, 2004.

\bibitem[{Mil}17]{milman}
E.~{Milman}.
\newblock {Beyond traditional curvature-dimension. I: New model spaces for
  isoperimetric and concentration inequalities in negative dimension.}
\newblock {\em {Trans. Am. Math. Soc.}}, 369(5):3605--3637, 2017.

\bibitem[MW82]{weissler}
C.~E. {Mueller} and F.~B. {Weissler}.
\newblock {Hypercontractivity for the heat semigroup for ultraspherical
  polynomials and on the n-sphere.}
\newblock {\em {J. Funct. Anal.}}, 48:252--283, 1982.

\bibitem[{Ngu}14]{nguyen-14}
V.~H. {Nguyen}.
\newblock {Dimensional variance inequalities of Brascamp-Lieb type and a local
  approach to dimensional Pr\'ekopa's theorem.}
\newblock {\em {J. Funct. Anal.}}, 266(2):931--955, 2014.

\bibitem[Ngu18]{18-nguyen}
V.~H. Nguyen.
\newblock Phi-entropy inequalities and asymmetric covariance estimates for
  convex measures.
\newblock 2018.
\newblock To appear in Bernoulli.

\bibitem[{Oht}16]{ohta}
S-{i} {Ohta}.
\newblock {$(K,N)$-convexity and the curvature-dimension condition for negative
  $N$.}
\newblock {\em {J. Geom. Anal.}}, 26(3):2067--2096, 2016.

\bibitem[{{\O}k}03]{oksendal}
B.~{{\O}ksendal}.
\newblock {\em {Stochastic differential equations. An introduction with
  applications.}}
\newblock Berlin: Springer, 6th ed. edition, 2003.

\bibitem[{Pav}14]{pavliotis}
G.~A. {Pavliotis}.
\newblock {\em {Stochastic processes and applications. Diffusion processes, the
  Fokker-Planck and Langevin equations.}}
\newblock New York, NY: Springer, 2014.

\bibitem[Sch01]{scheffer-phd}
G.~Scheffer.
\newblock {\em In\'egalit\'es fonctionnelles, g\'eom\'etrie conforme et noyaux
  markoviens}.
\newblock Phd, 2001.

\bibitem[{Sch}03]{scheffer}
G.~{Scheffer}.
\newblock {Local Poincar\'e inequalities in non-negative curvature and finite
  dimension.}
\newblock {\em {J. Funct. Anal.}}, 198(1):197--228, 2003.

\end{thebibliography}
\newcommand{\etalchar}[1]{$^{#1}$}

\end{document}